\documentclass[a4paper]{amsart}

\usepackage[english]{babel}
\usepackage{amsfonts}
\usepackage{amsmath}
\usepackage{amssymb}
\usepackage{amsthm}
\usepackage{latexsym}
\usepackage{amstext}
\usepackage{enumerate}
\usepackage{graphics}
\usepackage[all]{xy}
\usepackage{url}
\usepackage{hyperref}

\theoremstyle{plain}
\newtheorem{theorem}{Theorem}[section]

\newtheorem{lemma}[theorem]{Lemma}
\newtheorem{proposition}[theorem]{Proposition}

\theoremstyle{definition}
\newtheorem{definition}[theorem]{Definition}

\newtheorem{remark}[theorem]{Remark}
\numberwithin{equation}{section}

\newcommand{\N}{\mathbb{N}}

\newcommand{\Sf}{\mathbb{S}}
\newcommand{\T}{\mathbb{T}}

\begin{document}

\title[Principal torus bundles of Lorentzian $\mathcal{S}$-manifolds and \dots]
{Principal torus bundles of Lorentzian $\mathcal{S}$-manifolds and the $\varphi$-null Osserman condition}
\author{Letizia Brunetti \and Angelo V. Caldarella}

\address{Department of Mathematics, University of Bari ``Aldo Moro''\newline%
\indent Via E.~Orabona, 4 -- 70125 Bari (Italy)}%
\email[L.~Brunetti]{brunetti@dm.uniba.it}%
\email[A.V.~Caldarella]{caldarella@dm.uniba.it}%

\date{\today}

\begin{abstract}
The main result we give in this brief note relates, under suitable hypotheses, the $\varphi$-null Osserman, the null Osserman and the classical Osserman conditions to each other, via semi-Riemannian submersions as projection maps of principal torus bundles arising from a Lorentzian $\mathcal{S}$-manifold.
\end{abstract}

\maketitle

\section{Introduction}
The Jacobi operator is, for several reasons, one of the most interesting objects induced by the curvature operator.

On a (semi-)Riemannian manifold $(M,g)$, let us consider the unit spacelike $S^+(M)$ (resp.\ timelike $S^-(M)$) sphere bundle with fiber
\[
S^{\pm}_p(M)=\{z\in T_pM\,|\,g_p(z,z)=\pm1\},
\]
and put $S(M)=\bigcup_{p\in M}S_p^+(M)\cup S_p^-(M)$.

For any $z\in S_p(M)$, $p\in M$, the \emph{Jacobi operator with respect to $z$} is the endomorphism $R_z\colon z^{\bot}\rightarrow z^{\bot}$ such that $R_z(\cdot)=R_p(\cdot,z)z$ (\cite{GKV}), where $R$ is the $(1,3)$-type curvature tensor on $(M,g)$.

The Jacobi operator is obviously self-adjoint, hence a great deal of study has been carried out about the behaviour of its eigenvalues in the Riemannian case since R.\ Osserman proposed his Conjecture in \cite{OS} (see also \cite{Oss}).
Indeed, one easily sees that Riemannian space-forms are characterized by having Jacobi operators with exactly one constant eigenvalue corresponding to the sectional curvature. Those Riemannian manifolds whose Jacobi operators have eigenvalues independent both of the vector $z\in S_p(M)$ and of the point $p\in M$ are the \emph{Osserman manifolds}.
Any locally flat or locally rank-one symmetric space is an Osserman manifold, whilst the converse statement is known as the Osserman Conjecture. Several authors have dealt with this Conjecture, providing positive answers in many cases (\cite{C0},\cite{C1},\cite{C2},\cite{N01},\cite{N02},\cite{N03}).

One gets a different situation when considering the indefinite setting, where a semi-Riemannian manifold $(M,g)$ is said to be \emph{spacelike} (resp.\ \emph{timelike}) \emph{Osserman}, if the characteristic polynomial of $R_z$ is independent of both $z\in S^+_p(M)$ (resp.\ $z\in S^-_p(M)$) and $p\in M$. It is known that $(M,g)$ being spacelike Osserman is equivalent to $(M,g)$ being timelike Osserman (\cite{GKVV},\cite{GKV}), but several counterexamples to the Osserman Conjecture were found (see for example \cite{BBR}, \cite{BCGHV}, \cite{GVV}) for non-Lorentzian semi-Riemannian manifolds.

Finally, in the Lorentzian setting a complete solution for the Osserman Conjecture was provided in a sequence of works by E.~Garc\'ia-R\'io, D.N.~Kupeli and M.E.~V\'azquez-Abal (\cite{GK},\cite{GKVA}), together with N.~Bla\v{z}i\'c, N.~Bokan and P.~Gilkey (\cite{BBG}). They proved that a Lorentzian manifold is Osserman if and only if it has constant sectional curvature (see also \cite{GKV}).

It was defined a very fruitful new Osserman-related condition for Lorentzian manifolds in \cite{GKVA}. There, the authors introduced the Jacobi operator $\bar{R}_u$ with respect to a null (lightlike) vector $u$, and then they studied the so-called \emph{null Osserman conditions} with respect to a unit timelike vector (see also \cite{GKV}).

Here, we are concerned with an Osserman-related condition derived by the null Osserman condition, which is known as the \emph{$\varphi$-null Osserman condition}, introduced and studied by the first author in \cite{BR} for manifolds carrying Lorentzian globally framed $f$-structures. This condition appears to be a natural generalization of the null Osserman condition, to which it reduces when considering Lorentzian almost contact structures. This was motivated by the following considerations: although any Lorentzian Sasaki manifold $(M,\varphi,\xi,\eta,g)$ with constant $\varphi$-sectional curvature is globally null Osserman with respect to the timelike vector field $\xi$, there is no similar result when we consider Lorentzian $\mathcal{S}$-manifolds, which generalize Lorentzian Sasaki ones, and moreover, as we proved in \cite{BC}, no Lorentzian $\mathcal{S}$-manifold can be neither null Osserman, nor Osserman. Further basic properties of such manifolds are studied in \cite{BR} and developed in \cite{BC}.  We refer the reader to both works for more details about the $\varphi$-null Osserman condition, whilst the general reference for the whole Osserman framework is \cite{GKV}.

In this short note, we deal with the study of some relationships among the above three Osserman-related notions, providing a few results of equivalence, obtained by considering a natural structure of principal torus bundle arising from a Lorentzian $\mathcal{S}$-manifold, which involves semi-Riemannian submersions.

Indeed, from \cite{BLY}, a strong link between $f$-structures and Riemannian submersions is well-known. Namely, any compact and connected manifold endowed with a regular and normal $g.f.f$-structure is the total space of a torus principal bundle over a complex manifold, which, under suitable hypotheses, can be a K\"ahler manifold. Moreover, as also proved in \cite{BLY}, a compact, connected and regular Riemannian $\mathcal{S}$-manifold $(M,\varphi,\xi_\alpha,\eta^\alpha,g)$, with each $\xi_\alpha$ regular, projects itself onto a compact K\"ahler manifold and onto a compact and regular Sasakian manifold. These results have been extended to the semi-Riemannian case by the first author, together with A.M.~Pastore, who in \cite{BP} proved that a compact, connected and regular indefinite (in particular, Lorentzian) $\mathcal{S}$-manifold $(M,\varphi,\xi_\alpha,\eta^\alpha,g)$ projects itself onto a compact (indefinite) K\"ahler manifold and onto a compact and regular indefinite (Lorentzian) Sasakian manifold, via semi-Riemannian submersions.

Based on the above, after recalling, in Section $2$, some basic features of (almost) $\mathcal{S}$-manifolds, in Section $4$ we carry on an investigation on the possibilities of projectability of the $\varphi$-null Osserman conditions via semi-Riemannian submersions with a Lorentzian $\mathcal{S}$-manifold as total space, and either a Lorentzian Sasakian manifold or a K\"ahler manifold as base space. Using some properties established in \cite{BC}, which we briefly recall in Section $3$, together with a few properties of semi-Riemannian submersions, and under an additional assumption on the eigenvectors of the Jacobi operators, we obtain equivalence results relating the $\varphi$-null Osserman condition with the classical and the null Osserman condition in the framework of principal torus bundles constructed on a given Lorentzian $\mathcal{S}$-manifold.

In what follows, all smooth manifolds are supposed to be connected, and all tensor fields and maps are assumed to be smooth. Moreover, according to \cite{KN}, for the Riemannian curvature tensor of a semi-Riemannian manifold $(M,g)$ we use the definition $R(X,Y,Z,W)=g(R(Z,W)Y,X)=g(([\nabla_Z,\nabla_W]-\nabla_{[Z,W]})Y,X)$ for any vector fields $X,Y,Z,W$ on $M$.

Finally, for any $p\in M$ and any linearly independent vectors $x,y\in T_pM$ spanning a non-degenerate plane $\pi=\mathrm{span}(x,y)$, that is $g_p(x,x)g_p(y,y)-g_p(x,y)^2\neq0$, the sectional curvature of $(M,g)$ at $p$ with respect to $\pi$ is, by definition, the real number
\[
k_p(\pi)=k_p(x,y)=\frac{R_p(x,y,x,y)}{\Delta(\pi)},
\]
where $\Delta(\pi)=g_p(x,x)g_p(y,y)-g_p(x,y)^2$.

\section{Preliminaries}
Let us recall some basic definitions and facts about (almost) $\mathcal{S}$-manifolds needed in the rest of the paper.

Framed $f$-manifolds were originally considered by H.~Nakagawa in \cite{Na01} and \cite{Na02}, based on the notion of $f$-structure, which was firstly introduced in 1963 by K.~Yano (\cite{Y01}) as a generalization of both (almost) contact and (almost) complex structures. Such structures were later studied and developed by S.I.~Goldberg and K.~Yano (see, for example, \cite{GY}, \cite{GY02}) and, in the subsequent years, by several authors (\cite{Bl01}, \cite{BLY}, \cite{CFF}, \cite{IS01}, \cite{Van}).

A \emph{globally framed $f$-structure} (briefly \emph{$g.f.f$-structure}) on a manifold $M$ is a non-vanishing $(1,1)$-type tensor field $\varphi$ on $M$ of constant rank satisfying the following conditions: $\varphi^3+\varphi=0$, and the subbundle $\ker(\varphi)$ is parallelizable. This is equivalent to the existence of $s$ linearly independent vector fields $\xi_\alpha$ and $1$-forms $\eta^\alpha$ ($\alpha\in\{1,\ldots,s\}$), $s$ being the dimension of $\ker(\varphi))$ at any point $p\in M$, such that
\begin{equation}\label{eq:02}
\varphi^2=-I+\eta^\alpha\otimes\xi_\alpha\qquad\text{and}\qquad\eta^\alpha(\xi_\beta)=\delta^\alpha_\beta.
\end{equation}
Each $\xi_\alpha$ is said to be a \emph{characteristic vector field} of the structure, and a manifold $M$ carrying a $g.f.f$-structure is denoted by $(M,\varphi,\xi_\alpha,\eta^\alpha)$, and called a \emph{$g.f.f$-manifold}. When $s=1$ (resp.: $s=0$), we have an almost contact (resp.: almost complex) structure. From (\ref{eq:02}) one easily has $\varphi\xi_\alpha=0$ and $\eta^\alpha\circ\varphi=0$, for any $\alpha\in\{1,\ldots,s\}$. Furthermore, $\mathrm{Im}(\varphi)$ is a distribution on $M$ of even rank $r=2n$ on which $\varphi$ acts as an almost complex tensor field, and one has the splitting $TM=\mathrm{Im}(\varphi)\oplus\mathrm{ker}(\varphi)$, hence $\dim(M)=2n+s$. A $g.f.f$-manifold is said to be \emph{normal} if the $(1,2)$-type tensor field $N=[\varphi,\varphi]+2d\eta^\alpha\otimes\xi_\alpha$ vanishes.

In \cite{LP}, the authors study the properties of a $g.f.f$-manifold $(M,\varphi,\xi_\alpha,\eta^\alpha)$ endowed with a compatible indefinite metric, that is a semi-Riemannian metric $g$ verifying
\begin{equation}\label{eq:03}
	g(\varphi X,\varphi Y) = g(X,Y)-\sum_{\alpha=1}^s\varepsilon_\alpha\eta^\alpha(X)\eta^\alpha(Y),
\end{equation}
for all $X,Y\in\Gamma(TM)$, where $\varepsilon_\alpha=g(\xi_\alpha,\xi_\alpha)=\pm1$. Such a manifold is said to be an \emph{indefinite metric $g.f.f$-manifold} and denoted by $(M,\varphi,\xi_\alpha,\eta^\alpha,g)$. From (\ref{eq:03}) one also has $g(X,\xi_\alpha)=\varepsilon_\alpha\eta^\alpha(X)$ and $g(X,\varphi Y)=-g(\varphi X,Y)$, for any $X,Y\in\Gamma(TM)$, and the splitting $TM=\mathrm{Im}(\varphi)\oplus\mathrm{ker}(\varphi)$ becomes orthogonal.

The \emph{fundamental $2$-form} $\Phi$ of an indefinite metric $g.f.f$-manifold $(M,\varphi,\xi_\alpha,\eta^\alpha,g)$ is defined by $\Phi(X,Y)=g(X,\varphi Y)$. If $\Phi=d\eta^\alpha$, for any $\alpha\in\{1,\ldots,s\}$, the manifold is said to be an \emph{indefinite almost $\mathcal{S}$-manifold}. Finally, a normal indefinite almost $\mathcal{S}$-manifold is, by definition, an \emph{indefinite $\mathcal{S}$-manifold}. Such a manifold is characterized by the identity $(\nabla_X\varphi)Y=g(\varphi X,\varphi Y)\bar{\xi}+\bar{\eta}(Y)\varphi^2X$, where $\bar{\xi}=\sum_{\alpha=1}^s\xi_\alpha$ and $\bar{\eta}=\sum_{\alpha=1}^s\varepsilon_\alpha\eta^\alpha$. It follows that $\nabla_X\xi_\alpha=-\varepsilon_\alpha\varphi X$ and $\nabla_{\xi_\alpha}\xi_\beta=0$, for any $\alpha,\beta\in\{1,\ldots,s\}$, and each $\xi_\alpha$ is a Killing vector field.

For more details on (almost) $\mathcal{S}$-manifolds the reader is referred to \cite{DIP} in the Riemannian case, and to \cite{LP} for the indefinite case.

\section{Lorentzian $\mathcal{S}$-manifolds and the $\varphi$-null Osserman condition.}
The notion of $\varphi$-null Osserman condition is derived from that of null Osserman, which we briefly recall here, following \cite{GKVA} and \cite{GKV}.

Let $(M,g)$ be a Lorentzian manifold and $p\in M$. If $u\in T_pM$ is a lightlike (or null) vector, that is $u\neq0$ and $g_p(u,u)=0$, then $\mathrm{span}(u)\subset u^{\bot}$. We can endow the quotient space $\bar{u}^{\bot}=u^{\bot}/\mathrm{span}(u)$, whose canonical projection is $\pi\colon u^{\bot}\rightarrow\bar{u}^{\bot}$, with a positive definite inner product $\bar{g}$ defined by $\bar{g}(\bar{x},\bar{y})=g_p(x,y)$, where $\pi(x)=\bar{x}$ and $\pi(y)=\bar{y}$, obtaining the Euclidean vector space $(\bar{u}^{\bot},\bar{g})$.

The \emph{Jacobi operator with respect to $\bar u$} is the endomorphism $\bar R_u\colon\bar u^\bot \rightarrow\bar u^\bot$ defined by $\bar R_u(\bar x)=\pi(R_p(x,u)u)$, for all $\bar x=\pi(x)\in\bar u^\bot$. It is easy to see that $\bar R_u$ is a self-adjoint endomorphism, hence it is diagonalizable.

If $z\in T_pM$ is a unit timelike vector, the \emph{null congruence set} of $z$ is defined to be the set $N(z)=\{u\in T_pM\,|\,g_p(u,u)=0,\ g_p(u,z)=-1\}$. The elements of $N(z)$ are in one-to-one correspondence to those of the set $S(z)=\{x\in z^\bot\,|\,g_p(x,x)=1\}$, called the \emph{celestial sphere of $z$}, via the map $\psi\colon N(z)\rightarrow S(z)$ such that $\psi(u)=u-z$.

\begin{definition}[\cite{GKVA,GKV}]
A Lorentzian manifold $(M,g)$ is said to be \emph{null Osserman with respect to $z$}, $z\in T_pM$ being a unit timelike vector, if the eigenvalues of $\bar R_u$ and their multiplicities are independent of $u\in N(z)$.
\end{definition}

Following \cite{BR} and \cite{BC}, we recall the basic facts related with the definition of the $\varphi$-null Osserman condition.

Let $(M,\varphi,\xi_\alpha,\eta^\alpha,g)$ be a Lorentzian $g.f.f$-manifold, with $\dim(M)=2n+s$, and $\alpha\in\{1,\ldots,s\}$, $s\geqslant1$. It is easy to see that one of the characteristic vector fields has to be timelike and, without loss of generality, we assume it is $\xi_1$. If $p\in M$, we define the \emph{$\varphi$-celestial sphere} of $(\xi_1)_p$ to be the set $S_\varphi((\xi_1)_p)=S((\xi_1)_p)\cap\mathrm{Im}(\varphi_p)$, and the \emph{$\varphi$-null congruence set} of $(\xi_1)_p$ to be $N_{\varphi}((\xi_1)_p)=\psi^{-1}(S_\varphi((\xi_1)_p))$.

\begin{definition}[\cite{BR,BC}]
A Lorentzian $g.f.f$-manifold $(M,\varphi,\xi_\alpha,\eta^\alpha,g)$ is said to be \emph{$\varphi$-null Osserman with respect to $(\xi_1)_p$}, $p\in M$, if the eigenvalues of $\bar R_u$ and their multiplicities are independent of $u\in N_\varphi((\xi_1)_p)$.
\end{definition}

Fix $p\in M$ and consider $u\in N_\varphi((\xi_1)_p)$. Since we can write $u=(\xi_1)_p+x$, with $x\in S_\varphi((\xi_1)_p)$, there is a natural one-to-one correspondence between the two kinds of Jacobi operator $R_x\colon x^\bot\rightarrow x^\bot$ and $\bar R_u\colon\bar u^\bot\rightarrow\bar u^\bot$. In \cite{BC} it is provided the relationship between these two operators with respect to the $\varphi$-null Osserman condition, which we summarize in the following proposition.

\begin{proposition}[\cite{BC}]\label{prop:01}
Let $(M,\varphi,\xi_\alpha,\eta^\alpha,g)$ be a Lorentzian $\mathcal{S}$-manifold, $\dim(M)=2n+s$, $s\geqslant1$. For any $p\in M$, $M$ is $\varphi$-null Osserman with respect to $(\xi_1)_p$ if and only if the eigenvalues of $R_x$ with their multiplicities are independent of $x\in S_\varphi((\xi_1)_p)$.
\end{proposition}

The above result enables us to write the definition of the $\varphi$-null Osserman condition in terms of operator $R_x$, $x\in S_\varphi((\xi_1)_p)$, instead of $\bar R_u$, $u\in N_\varphi((\xi_1)_p)$. It is clear that, in the case of a Lorentzian Sasaki manifold, the $\varphi$-null Osserman condition reduces to that of null Osserman one.

\section{Principal torus bundles and the $\varphi$-null Osserman condition.}
From \cite{BLY} it is known that under an assumption of regularity it is possible to relate metric $g.f.f$-manifolds both to almost complex and to almost contact metric manifolds via Riemannian submersions. The semi-Riemannian version of the results of \cite{BLY} is provided in \cite{BP}, where it is possible to find the following result.

\begin{theorem}\label{teor:01}
Let $(M,\varphi,\xi_\alpha,\eta^\alpha,g)$ be a compact, connected and regular indefinite $\mathcal{S}$-manifold, with $\dim(M)=2n+s$, $s\geqslant2$.  Then, there exists a commutative diagram
\[
\xymatrix{
M\ar[dr]_\pi\ar[rr]^\tau&&M'\ar[dl]^{\pi'}\\
&N&
}
\]
where $N$ is a $2n$-dimensional compact K\"ahler manifold, either indefinite or not, and $M'$ is a $(2n+1)$-dimensional compact and regular Sasakian manifold, indefinite or not. All the maps are semi-Riemannian submersions with totally geodesic fibres, and more precisely:
\begin{itemize}
  \item $\tau$ is the projection of a principal $\T^{s-1}$-bundle over $M'$;
	\item $\pi'$ is the projection of a principal $\Sf^1$-bundle over $N$;
	\item $\pi$ is the projection of a principal $\T^s$-bundle over $N$.
\end{itemize}
where $\T^{k}$ is the $k$-dimensional torus, for any $k\in\N$, $k\geqslant1$.
\end{theorem}

For the notion of regularity of a distribution and of a $g.f.f$-structure the reader is referred to \cite{Pa} and \cite{BLY}. The general idea of this result, as contained in \cite{BLY}, is to fibrate $M$ by any $s-r$ of the vector fields $\xi_\alpha$'s, to obtain a principal $\T^{s-r}$-bundle over a $(2n+r)$-dimensional manifold $M'$. The remaining $r$ characteristic vector fields are then projectable to $M'$, inducing a $g.f.f$-structure on $M'$ and preserving the regularity. Thus, $M'$ can be fibrated again by its $r$ characteristic vector fields, obtaining a principal $\T^r$-bundle over $N$, which finally produces a commutative diagram. In particular, if we fibrate a Lorentzian $\mathcal{S}$-manifold $M$ by the $s-1$ spacelike characteristic vector fields, in Theorem \ref{teor:01} we obtain that $N$ is a K\"ahler manifold and $M'$ is a Lorentz Sasakian manifold.

We are going to find out some informations about the possibility of projecting the $\varphi$-null Osserman condition both onto the null Osserman condition and the classical Osserman condition, via the previous fibrations.

In general (see \cite{FIP}, \cite{ONArt}), given a $C^{\infty}$-submersion $f:(M,g)\rightarrow(B,g')$ between semi-Riemannian manifolds, i.e.\ a map whose differential $(df)_p$ is surjective, for all $p\in M$, then $\mathcal{V}=(\ker(df)_p)_{p\in M}$ and $\mathcal{H}=(\ker(df)_p^\bot)_{p\in M}$ are, by definition, the \emph{vertical} and the \emph{horizontal} distributions of $f$. Such a map is said to be a \emph{semi-Riemannian submersion} if each fibre $f^{-1}(p')$, $p'\in B$, is a (semi-)Riemannian submanifold of $M$ and the restriction of $g_p$ to $\mathcal{H}_p$ is an isometry for all $p\in M$. A vector field $U$ (resp.\ $X$) on $M$ such that $U_p\in\mathcal{V}_p$ (resp.\ $X_p\in\mathcal{H}_p$) is called \emph{vertical} (resp.\ \emph{horizontal}). A vector field $X$ on $M$ such that there exists a vector field $X'$ on $B$ for which $f_*X=X'$ is said to be \emph{projectable}, and any horizontal, projectable vector field on $M$ is said to be \emph{basic}. The vertical distribution is always integrable, with the fibres of $f$ as leaves. Denoting by $v$ and $h$ the projections of $TM$ onto $\mathcal{V}$ and $\mathcal{H}$, respectively, the \emph{O'Neill tensors} of $f$ are the $(1,2)$-type tensor fields $T$ and $A$ on $M$ defined by:
\[
\begin{split}
&T(X,Y)=T_XY:=v\nabla_{vX}hY+h\nabla_{vX}vY,\\
&A(X,Y)=A_XY:=v\nabla_{hX}hY+h\nabla_{hX}vY.
\end{split}
\]
They are both $g$-skew-symmmetric tensors, and they satisfy the following fundamental properties:
\[
\begin{array}{lll}
T_UW=T_WU&\qquad& U,W\in\mathcal{V}\\
A_XY=-A_YX=\frac{1}{2}v[X,Y]&\qquad& X,Y\in\mathcal{H}
\end{array}
\]
It follows that the horizontal distribution is integrable if and only if $A=0$, and in this case the leaves are totally geodesic submanifolds of $M$. Furthermore, the fibres of $f$ are totally geodesic semi-Riemannian submanifolds of $M$ if and only if $T=0$.

\begin{lemma}\label{lem:01}
Let $(M,\varphi,\xi_\alpha,\eta^\alpha,g)$ be a Lorentzian $\mathcal{S}$-manifold, with $\dim(M)=2n+s$, $s\geqslant1$. Let $\pi:M\rightarrow N$ be a principal $\T^s$-bundle over a K\"ahler manifold, as in Theorem \ref{teor:01}. We have:
\begin{equation}\label{eq:11}
A_XY=-g(X,\varphi Y)\bar\xi,\qquad A_X\xi_\alpha=-\varepsilon_\alpha\varphi X,
\end{equation}
for any $X,Y\in\mathrm{Im}(\varphi)$ and any $\alpha\in\{1,\ldots,s\}$, where $\bar\xi=\sum_{\alpha=1}^s\xi_\alpha$.
\end{lemma}
\proof By construction of $\pi$, we have $\mathcal{H}_p=\mathrm{Im}(\varphi_p)$ and $\mathcal{V}_p=\mathrm{span}((\xi_1)_p,\ldots,(\xi_s)_p)$ for any $p\in M$. Thus, since $\nabla_X\xi_\alpha=-\varepsilon_\alpha\varphi X$, by direct calculation we get:
\[
\begin{split}
A_XY&=
v(\nabla_XY)=
\sum_{\alpha=1}^s\varepsilon_\alpha g(\nabla_XY,\xi_\alpha)\xi_\alpha=
-\sum_{\alpha=1}^s\varepsilon_\alpha g(Y,\nabla_X\xi_\alpha)\xi_\alpha\\
&=\sum_{\alpha=1}^s g(Y,\varphi X)\xi_\alpha=
-g(X,\varphi Y)\bar\xi,
\end{split}
\]
for all $X,Y\in\mathcal{H}$. Analogously, we have $A_X\xi_\alpha=h(\nabla_X\xi_\alpha)=-\varepsilon_\alpha\varphi X$ for all $X\in\mathcal{H}$ and $\alpha\in\{1,\ldots,s\}$. \endproof

For a semi-Riemannian submersion $f\colon (M,g)\rightarrow (B,g')$, let us denote by $R^*$ the $(1,3)$-type $\mathcal{H}$-valued tensor field on $M$ such that, if $X,Y,Z\in\Gamma(TM)$ are basic vector fields $f$-related to $X',Y',Z'\in\Gamma(TB)$, then $R^*(X,Y)Z $ is the unique basic vector field $f$-related to $R'(X',Y')Z'$. Thus, for any $x\in\mathcal{H}_p$, one can consider the self-adjoint endomorphism $R^*_x\colon x^\bot\cap\mathcal{H}_p\rightarrow x^\bot\cap\mathcal{H}_p$ such that $R^*_x(y)=R^*_p(y,x)x$.

\begin{lemma}
Let $f\colon (M,g)\rightarrow (B,g')$ be a semi-Riemannian submersion. For any orthogonal vectors $x,y\in\mathcal{H}_p$ one has
\begin{equation}\label{eq:12}
R^*_x(y)=h_pR_x(y)-3A_xA_x(y).
\end{equation}
\end{lemma}
\proof From standard formulas on the curvature tensors of a submersion (see \cite{FIP}, pag.~13), we have
\[
\begin{split}
g_p(R^*_x(y),z)&=R^*_p(x,y,x,z)\\
               &=R_p(x,y,x,z)+2g_p(A_x(y),A_x(z))-g_p(A_y(x),A_x(z))\\
               &=g_p(h_pR_x(y),z)-3g_p(A_xA_x(y),z)
\end{split}
\]
for any $z\in x^\bot\cap\mathcal{H}_p$, which yields (\ref{eq:12}). \endproof

\begin{proposition}\label{prop:02}
Let $(M,\varphi,\xi_\alpha,\eta^\alpha,g)$ be a Lorentzian $\mathcal{S}$-manifold, with $\dim(M)=2n+s$, $s\geqslant1$. Let $\pi:M\rightarrow N$ be a principal $\T^s$-bundle over a K\"ahler manifold $(N,J,G)$, as in Theorem \ref{teor:01}. Let $p\in M$, and suppose that, for any $x\in S_\varphi((\xi_1)_p)$, $\varphi x$ is an eigenvector of $R_x$. Then, $M$ is $\varphi$-null Osserman with respect to $(\xi_1)_p$ if and only if $N$ is Osserman at $p'=\pi(p)$.
\end{proposition}
\proof Suppose first that $s\geqslant 2$. Fix $p'\in N$, with $p'=\pi(p)$, $p\in M$, and let $x'\in T_{p'}N$ a unit vector, and $y',z'\in x'^{\bot}$. Let $x\in S_\varphi((\xi_1)_p)$, $V=x^\bot\cap\mathrm{Im}(\varphi_p)$ and $y,z\in V$ such that $x'=(d\pi)_p(x)$, $y'=(d\pi)_p(y)$ and $z'=(d\pi)_p(z)$. Then
\[
g_p(R^*_x(y),z)=G_{p'}((d\pi)_p(R^*_x(y)),(d\pi)_p(z))=G_{p'}(R'_{x'}(y'),z'),
\]
which implies that the Jacobi operators $R^*_x\colon V\rightarrow V$ and $R'_{x'}\colon x'^{\bot}\rightarrow x'^{\bot}$ have the same characteristic polynomial. Using (\ref{eq:11}) one has $A_xA_x(y)=-(s-2)g_p(y,\varphi x)\varphi x$ and since $R_x$ leaves the subspace $V$ invariant, (\ref{eq:12}) gives
\[
R^*_x(y)=R_x(y)+3(s-2)g_p(y,\varphi x)\varphi x
\]
for any $y\in V$. Observe that if $\varphi x$ is an eigenvector of $R_x$, we have
\[
g_p(R_x(y),\varphi x)\varphi x=g_p(y,R_x(\varphi x))\varphi x=R_x(g_p(y,\varphi x)\varphi x),
\]
that is the endomorphism of $V$ such that $y\mapsto g_p(y,\varphi x)\varphi x$ commutes with $R_x$. This implies they are simultaneously diagonalizable, and if $\lambda_i$, $i\in\{1,\ldots,r\}$ are the eigenvalues of $R_x$, counted with multiplicities, with $\lambda_1$ relative to $\varphi x$, then $\lambda_1+3(s-2),\lambda_j$, $j\in\{2,\ldots,r\}$ are the eigenvalues of $R^*_x$. By Proposition \ref{prop:01} we obtain our statement.

If $s=1$ then the proof goes through as above, except for the fact that one has $A_xA_x(y)=g_p(y,\varphi x)\varphi x$. \endproof

\begin{remark}
It is clear, from the previous proof, that in case $s=2$ the hypothesis of $\varphi x$ being an eigenvector of $R_x$ can be dropped without affecting the result. Furthermore, in case $s=1$, the statement is relative to the projection $\pi'$ of the commutative diagram in the Theorem \ref{teor:01}.
\end{remark}

\begin{lemma}
Let $(M,\varphi,\xi_\alpha,\eta^\alpha,g)$ be a Lorentzian $\mathcal{S}$-manifold, with $\dim(M)=2n+s$, $s\geqslant2$. Let $\tau\colon M\rightarrow M'$ be a principal $\T^{s-1}$-bundle over a Lorentz Sasakian manifold, as in Theorem \ref{teor:01}. We have:
\begin{equation}\label{eq:15}
A_XY=-g(X,\varphi Y)\sum_{\alpha=2}^s\xi_\alpha,\qquad A_X\xi_\alpha=-\varphi X,
\end{equation}
for any $X,Y\in\mathrm{Im}(\varphi)\oplus\mathrm{span}(\xi_1)$ and any $\alpha\in\{2,\ldots,s\}$.
\end{lemma}
\proof By construction of $\tau$, we have the splitting $\mathcal{H}_p=\mathrm{Im}(\varphi_p)\oplus\mathrm{span}(\xi_1)$ and $\mathcal{V}_p=\mathrm{span}((\xi_2)_p,\ldots,(\xi_s)_p)$ for any $p\in M$. Proceeding along the same lines as the proof of Lemma \ref{lem:01}, we get \eqref{eq:15}. \endproof

\begin{proposition}\label{prop:03}
Let $(M,\varphi,\xi_\alpha,\eta^\alpha,g)$ be a Lorentzian $\mathcal{S}$-manifold, with $\dim(M)=2n+s$, $s\geqslant2$. Let $\tau:M\rightarrow M'$ be a principal $\T^{s-1}$-bundle over a Lorentz Sasakian manifold $M'$ with structure $(\varphi',\xi',\eta',g')$ as in Theorem \ref{teor:01}. Let $p\in M$, and suppose that, for any $x\in S_\varphi((\xi_1)_p)$, $\varphi x$ is an eigenvector of $R_x$. Then, $M$ is $\varphi$-null Osserman with respect to $(\xi_1)_p$ if and only if $M'$ is null Osserman with respect to $\xi'_{p'}$, $p'=\tau(p)$.
\end{proposition}
\proof One can follow the same proof of Proposition \ref{prop:02} where, using (\ref{eq:15}), one has $A_xA_x(y)=-(s-1)g_p(y,\varphi x)\varphi x$. \endproof

Propositions \ref{prop:02} and \ref{prop:03} can be summarized as follows.

\begin{theorem}
Let $(M,\varphi,\xi_\alpha,\eta^\alpha,g)$ be a compact, connected and regular Lorentzian $\mathcal{S}$-manifold, with $\dim(M)=2n+s$, $s\geqslant2$. Consider the commutative diagram of principal torus bundles
\[
\xymatrix{
M\ar[dr]_\pi\ar[rr]^\tau&&M'\ar[dl]^{\pi'}\\
&N&
}
\]
where $N$ is a $2n$-dimensional compact K\"ahler manifold and $M'$ is a $(2n+1)$-dimensional compact and regular Lorentz Sasakian manifold, with unit timelike characteristic vector field $\xi'=\tau_*(\xi_1)$. Let $p\in M$, and suppose that $\varphi x$ is an eigenvector of $R_x$ for any $x\in S_\varphi((\xi_1)_p)$. The following three statements are equivalent.
\begin{itemize}
  \item[\textrm{(a)}] $M$ is $\varphi$-null Osserman with respect to $(\xi_1)_p$;
	\item[\textrm{(b)}] $N$ is Osserman at $q=\pi(p)$;
	\item[\textrm{(c)}] $M'$ is null Osserman with respect to $\xi'_{p'}$, $p'=\tau(p)$.
\end{itemize}
\end{theorem}

\begin{remark}
It is clear that the three Osserman-type conditions in the above theorem can be also considered either pointwise or globally. Moreover, if we use the pointwise conditions, from the equivalence $(a)\Leftrightarrow(b)$ it follows that $N$ is Einstein at each point and the connectedness implies that it is a K\"ahler-Einstein manifold.
\end{remark}

\begin{remark}
In case $\tau:M\rightarrow M'$ is a principal $\T^{s-1}$-bundle from a Lorentzian $\mathcal{S}$-manifold $(M,\varphi,\xi_\alpha,\eta^\alpha,g)$ with $\dim(M)=2n+s$, $s\geqslant2$, over a Sasakian manifold $M'$ with structure $(\varphi',\xi',\eta',g')$ as in Theorem \ref{teor:01}, we could ask about the Osserman condition on $M'$. Let us suppose $M'$ pointwise Osserman, since it is odd-dimensional, it has constant sectional curvature $c$ (\cite{C0,GKV}). Being $k(X',\xi')=1$, for any $X'\in\mathrm{Im(\varphi')}$, then $c=1$ and $M'$ is locally isometric to the sphere $\Sf^{2n+1}$ with its standard Sasakian structure (see \cite{Bl1}, p.~114). By construction of the bundle projection $\tau$, we can suppose that $\mathcal{H}_p=\mathrm{Im}(\varphi_p)\oplus\mathrm{span}((\xi_s)_p)$ and $\mathcal{V}_p=\mathrm{span}((\xi_1)_p,\ldots,(\xi_{s-1})_p)$. Hence, with calculations similar to those of Lemma \ref{lem:01}, one has $A_XY=g(Y,\varphi X)\sum_{\alpha=1}^{s-1}\xi_\alpha$. By standard formulas on sectional curvatures of the total and the base spaces of a semi-Riemannian submersion (see \cite{FIP}, p.~14) we have
\[
k(x,\varphi x)=k'(x',\varphi' x')-3g(A_x\varphi x,A_x\varphi x)=1-3(s-3),\qquad x\in\mathrm{Im}(\varphi_p),
\]
which gives a necessary condition on the $\varphi$-sectional curvature of $M$ for $M'$ to be an Osserman Sasakian manifold. 
\end{remark}

\begin{remark}
Analogously, in case $\tau:M\rightarrow M'$ is a principal $\T^{s-1}$-bundle from a Lorentzian $\mathcal{S}$-manifold $(M,\varphi,\xi_\alpha,\eta^\alpha,g)$, with $\dim(M)=2n+s$, $s\geqslant2$, over a Lorentz Sasakian manifold $M'$ with structure $(\varphi',\xi',\eta',g')$ as in Theorem \ref{teor:01}, we could ask again about the Osserman condition on $M'$. It is known that any connected Lorentzian Osserman manifold is a space-form (\cite{GKV}), and since $k(X',\xi')=-1$, for any $X'\in\mathrm{Im(\varphi')}$, $M'$ has constant sectional curvature $c=-1$. As in the previous calculations, using \eqref{eq:15}, we have
\[
k(x,\varphi x)=k'(x',\varphi' x')-3g(A_x\varphi x,A_x\varphi x)=-1-3(s-1),\qquad x\in\mathrm{Im}(\varphi_p),
\]
which is a necessary condition on the $\varphi$-sectional curvature of $M$ for $M'$ to be a Lorentzian Osserman manifold.
\end{remark}

\end{document}